\documentclass[12pt]{article}
\usepackage[a4paper,margin=28mm]{geometry}
\usepackage{amsmath,amssymb,amsthm,mathtools,microtype,enumitem,booktabs}
\PassOptionsToPackage{hyphens}{url}
\usepackage[hidelinks]{hyperref}
\hypersetup{pdftitle={One Unit Separates Polynomial Time from Undecidability in Term Coding},pdfauthor={Soren Riis}}
\newtheorem{theorem}{Theorem}[section]
\newtheorem{lemma}[theorem]{Lemma}
\newtheorem{proposition}[theorem]{Proposition}
\newtheorem{corollary}[theorem]{Corollary}
\theoremstyle{definition}
\newtheorem{definition}[theorem]{Definition}
\newtheorem{example}[theorem]{Example}

\newcommand{\DD}{D}

\newcommand{\Sym}{\operatorname{Sym}}
\newcommand{\id}{\operatorname{id}}

\newcommand{\PaperThreeLeanGitHub}{\url{https://github.com/SR123/term-coding-disequality-lean}}
\newcommand{\PaperThreeLeanDOI}{\href{https://doi.org/10.5281/zenodo.22727895}{\texttt{10.5281/zenodo.22727895}}}
\newcommand{\PaperThreeInteractive}{\url{https://sr123.github.io/term-coding-disequality-lean/}}
\title{One Unit Separates Polynomial Time\\from Undecidability in Term Coding}
\author{S\o ren Riis}
\date{12 September 2026}
\begin{document}
\maketitle
\begin{abstract}
Term coding provides a common algebraic framework for network coding,
index coding and problems in extremal combinatorics. We exhibit a
decision problem in which lowering an output threshold by just one
changes the complexity from polynomial time to undecidability: no
algorithm then halts with the correct answer on every input.
The problem concerns dispersion, the maximum number of distinct output
tuples obtainable by interpreting the function symbols in a tuple of
terms on a finite alphabet. We restrict inputs by inequalities between
terms and ask whether a given instance meets a prescribed output
threshold for some alphabet size at least two. Both thresholds are
considered on the same class of instances.
A machine-checked Lean development and an interactive presentation of
the paper accompany this work
(GitHub: \PaperThreeLeanGitHub; DOI: \PaperThreeLeanDOI);
Section~\ref{sec:open} specifies its external input and verification limits.
\end{abstract}

\section{Two decision problems with adjacent thresholds}
An expression such as $f(g(x),y)$ becomes a function of $x,y$ once a set
$A$ and functions $g:A\to A$ and $f:A^2\to A$ are chosen. We call such
expressions \emph{terms}, and the choice of functions an
\emph{interpretation}. Every occurrence of a symbol uses the same
function: for example, $g(g(x))$ means that this function is applied
twice. A tuple of terms in variables $x_1,\ldots,x_m$ defines a map
from $A^m$ to $A^r$, where $r$ is the number of terms in the tuple.
We restrict its domain by requiring finitely many inequalities
$u\ne v$ between terms. Such inequalities are also called
\emph{disequalities} or \emph{tests}.

An instance $\Gamma$ specifies the variables, function symbols and their
arities, output terms, and tests. For each $n\ge1$, let $D_\Gamma(n)$
be the maximum number of distinct output tuples obtained from inputs
satisfying all tests, as the interpreting functions on
$[n]=\{0,\ldots,n-1\}$ vary. This maximum is called the
\emph{dispersion}. It counts output tuples, not satisfying input tuples:
different inputs can give the same output. The finite set on which the
functions act is also called the \emph{alphabet}. The instance specifies
neither its size nor any of the interpreting functions.
Section~\ref{sec:language} gives the formal definitions.

Let $\mathcal C$ consist of all instances with a finite set $X$ of
variables containing distinct variables $x,y$, output tuple
\[
 (x,y,t_3,\ldots,t_r),\qquad r\ge2,
\]
and a finite list of tests that includes $x\ne y$. The remaining output
terms may use any variables in $X$ and may repeat. All variables, symbols,
terms and tests are explicitly listed in the input. The class
$\mathcal C$ does not depend on the integer $k$ used below; in particular,
the numbers of variables and outputs need not equal $k$.

For $k\ge2$, put
\[
 b_k(n)=n^k-n^{k-1}=n^{k-1}(n-1).
\]
This is the number of $k$-tuples in $[n]^k$ whose first two coordinates
differ. We compare two questions about an input $\Gamma\in\mathcal C$:
does $D_\Gamma(n)\ge b_k(n)$ hold for some $n\ge2$, and does
$D_\Gamma(n)\ge b_k(n)+1$ hold for some $n\ge2$? Both quantify over
the same possible sizes of the finite set. The second question asks for
just one more output tuple, but admits a polynomial-time algorithm,
whereas the first is undecidable.

The algorithm uses a graph determined by the terms. Its vertices are
the distinct subterms, together with all declared variables; an edge
runs from an immediate argument to the term containing it. This is a
directed acyclic graph. A \emph{term cut} is a set of vertices meeting
every directed path from a variable to an output term. Paths of length
zero are included, and a cut may contain variables or output terms.
Let $\rho_\Gamma$ be the minimum size of such a cut. Every cut for an
instance in $\mathcal C$ contains $x,y$, since each is both a variable
and an output. Thus $\rho_\Gamma\ge2$.

\begin{theorem}[The one-unit difference in complexity]\label{thm:main}
For every fixed integer $k\ge2$, the following decision problems on
$\mathcal C$ have the indicated complexities:
\[
\begin{array}{lll}
 \text{input: }\Gamma\in\mathcal C&
 \exists n\ge2:\ D_\Gamma(n)\ge b_k(n)&\text{undecidable},\\[2pt]
 \text{input: }\Gamma\in\mathcal C&
 \exists n\ge2:\ D_\Gamma(n)\ge b_k(n)+1&\text{polynomial time}.
\end{array}
\]
The second answer is yes if and only if no test has syntactically
identical sides and $\rho_\Gamma\ge k+1$. This criterion can also be
decided in polynomial time when $k\ge2$ is part of the input, written
in binary. If it holds, and $s$ is the number of distinct subterms in
the outputs and tests together with all declared variables, then
\[
 D_\Gamma(n)\ge\lfloor n/s\rfloor^{k+1}\quad(n\ge s).
\]
In particular, whenever the criterion holds, the strict threshold is
attained at the alphabet size $n_0=s^{k+1}$.
For every fixed $k$, the first problem remains undecidable on the
subclass with $\rho_\Gamma=k$, where $D_\Gamma(n)\le b_k(n)$ for all
$n\ge2$. There is no total computable function $B$ such that every
positive instance in this subclass, of input length $N$, has a size
$n$ with $2\le n\le B(N)$ and $D_\Gamma(n)=b_k(n)$.
\end{theorem}

We call the two inequalities in the theorem the \emph{non-strict} and
\emph{strict} thresholds, respectively: since $D_\Gamma(n)$ is an
integer, the second is equivalent to $D_\Gamma(n)>b_k(n)$.
For fixed $k$, the sufficient size $s^{k+1}$ is polynomial in $s$.
When $k$ varies, its value need not be polynomial, but its binary
representation is still short. On a positive strict instance,
$k+1\le\rho_\Gamma\le s$, so it has at most
$s(\lfloor\log_2 s\rfloor+1)$ bits. This bounds the length of the
integer specifying a sufficient alphabet size; it does not assert that
all entries of functions of arbitrary arity can be listed in polynomial
time.

\paragraph{Why all output pairs force identities.}
For outputs $(x,y)$, the test $x\ne y$ allows at most $n(n-1)$ outputs.
Achieving this number means that every pair $(a,b)$ with $a\ne b$
satisfies every test. In particular, if a further test is $F(x)\ne y$,
then $F(a)=a$ for every $a$: otherwise the pair $(a,F(a))$ would fail.
We refer to equality with this upper bound as \emph{saturation}.
Applying this observation to terms expressing inverse and relator laws
forces them to hold on the whole finite set. Section~\ref{sec:reduction}
uses it to encode the following undecidable problem for one fixed
finitely presented group: given a word, does it have nonidentity image
in some finite quotient? The external group-theoretic result is
\cite[Theorem 2.1]{BW}. Only the term representing the input word varies;
the unary symbols and the number of tests are fixed.
Adding $k-2$ new variables as further outputs, with no tests involving
them, multiplies the image size by $n^{k-2}$ and proves the other degrees.

\paragraph{The cut bound and earlier routing work.}
Riis and Gadouleau \cite{RG} proved that ordinary dispersion, without
tests, has growth exponent equal to the minimum term cut. Their lower
bound transmits independent values along disjoint paths in the term
graph. Each transmitted value is accompanied by a label identifying
its subterm, so that different uses of a symbol can be handled by one
function. This construction is called \emph{routing}.
Section~\ref{sec:routing} shows that the same labels satisfy every test
whose sides are syntactically different. It also proves the sharper
upper bound
\[
 D_\Gamma(n)\le n^{\rho_\Gamma}-n^{\rho_\Gamma-1}
 \qquad(\Gamma\in\mathcal C),
\]
because every cut contains the output variables $x,y$, whose values
must differ. These two bounds give the strict decision criterion.
The reduction uses the new identity-forcing argument above; it does
not follow from the earlier asymptotic theorem.

\paragraph{Approaching the upper bound by explicit functions.}
For the two-variable instance associated to a fixed presentation $G$
and a word $w$, write $D_{G,w}(n)$ for the dispersion.
For every nonempty word $w$, Section~\ref{sec:near} constructs functions
satisfying
\[
 n(n-1)-Mn\le D_{G,w}(n)\le n(n-1)\qquad(n>|w|),
\]
where $|w|$ counts the letters as written and $M$ depends only on the
presentation. Given $w$ and a positive rational $\varepsilon$, the
construction computes an integer $n\ge2$ and unary functions for which
the allowed output image contains at least $(1-\varepsilon)n(n-1)$ pairs. It works even
for words for which equality with $n(n-1)$ is impossible at every size.
Thus arbitrarily accurate proportions of this upper bound, with
computable witnesses, do not decide whether the bound can be reached.
The reason is concrete: each failed identity at a point excludes at
most one further second coordinate. At most $Mn$ pairs are lost, a
vanishing proportion of $n(n-1)$, whereas losing no pair forces every
identity to hold everywhere. The approximate constructions need not
satisfy the group laws.

\begin{example}[Both answers occur at every degree]\label{ex:nonvacuity}
Fix $k\ge2$, variables $x_1,\ldots,x_{k+1}$, and only the test
$x_1\ne x_2$. The output tuple $(x_1,\ldots,x_{k+1})$ has
$b_{k+1}(n)$ distinct outputs and exceeds $b_k(n)+1$ already at $n=2$.
With the same variables, test and number of output coordinates, the
tuple $(x_1,\ldots,x_k,x_1)$ instead has exactly $b_k(n)$ outputs.
It gives a no answer to the strict question and a yes answer to the
non-strict question. Adding the deliberately impossible test $x_1\ne x_1$
excludes every assignment and gives $D_\Gamma(n)=0$, so both questions
have answer ``no''.
For $k=2$, the first two examples are $(x,y,z)$ with
$D(n)=n^2(n-1)$ and $(x,y,x)$ with $D(n)=n(n-1)$.
\end{example}

The strict positive instances form a subset of the non-strict positive
instances; there is no contradiction in the subset being easier to
recognise. Also $b_k(n)+1<b_{k+1}(n)$ for every $n\ge2$.
As $k$ increases, the undecidable thresholds therefore alternate with
the polynomial-time thresholds. The input class remains $\mathcal C$
throughout.

\paragraph{Correction to the earlier omnibus.}
This paper is intended to replace \emph{Term Coding for Extremal
Combinatorics: Dispersion and Complexity Dichotomies}, arXiv:2504.16265
\cite{Omnibus}. That work has been divided into focused papers.
The entropic framework and guessing-number sandwich theorem are developed
in \cite{Paper1}; ordinary dispersion and its exact and asymptotic
decision problems are studied in \cite{Paper2}.
Versions~1--3 of the earlier omnibus asserted undecidability of finite
square term bijectivity, the existence of an interpretation making a
tuple of $k$ terms in $k$ variables a bijection on $A^k$ for some finite
$A$ with $|A|\ge2$. The proposed Horn-clause reduction in Appendix~A
of v1 did not establish the required forcing properties, and the
machine-simulation sketches in Section~7.1 of v2 and v3 did not establish
the required reduction. The undecidability assertion and the unconditional
pure-equation complexity conclusions depending on it are therefore
withdrawn. Finite square term bijectivity remains open. The present
reduction applies to the stated class with disequality tests and does
not establish that earlier assertion.

\section{Terms, interpretations, and the output image}\label{sec:language}
Fix a finite nonempty set $X$ of variables, also called \emph{source
variables}, and a finite signature $\Sigma$: a set of function symbols
with a specified nonnegative integer arity for each symbol. Terms are
finite rooted expressions over $X$ and $\Sigma$, with ordered arguments.
A symbol of arity zero is a constant symbol. Syntactic equality means
literal equality of expressions using the same symbol and variable
names; no algebraic identities or cancellation rules are imposed.
An interpretation $I$ on a nonempty set $A$ assigns one function
$I_f:A^{\operatorname{arity}(f)}\to A$ to each $f\in\Sigma$.
Every such function is defined on its entire domain, with no prescribed
values. For an assignment $a:X\to A$, the value $u^I(a)$ of a term is
defined recursively. All occurrences of $f$ use the same $I_f$.

\begin{definition}
An instance $\Gamma$ consists of variables $X=\{x_1,\ldots,x_m\}$,
a finite signature $\Sigma$, an output tuple
$\mathbf t=(t_1,\ldots,t_r)$, and a finite list $\Delta$ of tests
$u\ne v$ between terms over $X$ and $\Sigma$.
For an interpretation $I$ on a finite nonempty set $A$, define
\[
 C_I=\{a\in A^m:u^I(a)\ne v^I(a)
                    \text{ for every }(u\ne v)\in\Delta\}.
\]
The output map is
$\mathbf t^I(a)=(t_1^I(a),\ldots,t_r^I(a))$.
Its image on the allowed inputs is $\mathbf t^I(C_I)\subseteq A^r$.
For $n\ge1$, the dispersion is
\[
 D_\Gamma(n)=\DD_n(\mathbf t,\Delta)
       =\max_{I\text{ on }[n]}|\mathbf t^I(C_I)|.
\]
The notation $\DD_n(\mathbf t,\Delta)$ leaves $X$ and $\Sigma$ implicit.
\end{definition}

A bijection between two finite sets transports interpretations and
assignments, preserving the number of outputs. Thus only $|A|$ matters.
There are finitely many functions of each specified arity on $[n]$,
so the maximum exists. In particular, $D_\Gamma(n)\ge h$ means that
some interpretation on $[n]$ has at least $h$ distinct outputs.

For outputs $(x,y)$ and only the test $x\ne y$, the output map is the
identity on pairs and has $n(n-1)$ outputs. Input and output counts can
differ when other variables are present. For example, with variables
$x,y,z$, outputs $(x,y,x)$ and only the test $x\ne y$, there are
$n^2(n-1)$ allowed inputs, but only $n(n-1)$ different outputs:
changing $z$ does not change the output.

\paragraph{Writing nested terms as acyclic equations.}
One may introduce a new variable for each distinct nonvariable subterm
and an equation recording its computation. For example, the output
$f(g(x),y)$ can be written using $z=g(x)$ and $v=f(z,y)$, with output
$v$. Every assignment to the original variables has a unique extension
to the new ones. Replacing terms in outputs and tests by their names
therefore preserves the output image. Equal subterms receive the same
name, and every use of a symbol still has the same interpretation.
These equations merely record evaluation; they add no further equality
condition on the original variables. This rewriting is often called
\emph{flattening}.

\paragraph{Finite input descriptions.}
For the complexity statements, names, arities, variables, output and
test lists, and ordered argument lists are explicitly encoded in binary.
Integers have canonical self-delimiting encodings, and every list includes
its length and all entries. Let $N$ be the bit length of this description.
Terms may be represented by syntax trees or by an acyclic graph listed
in topological order, with arguments referring to earlier nodes.
In the latter representation, different nodes can still represent the
same expression. We identify them by syntax, processing child terms
first, without expanding the graph into a potentially exponential tree.
Validation and this identification take polynomial time in either
representation.

The algorithm rejects invalid references, wrong argument counts,
noncanonical integer encodings and trailing data. A declared list length
or arity does not justify a long search for absent entries: a parser
stops and rejects when the supplied data end. Thus validation is
polynomial in the actual input length even for malformed descriptions.
For the uniform problem, $k$ is a separate binary integer and inputs
with $k<2$ are rejected.

\section{Image bounds from cuts in the term graph}\label{sec:routing}
Let $S$ contain all distinct subterms of the outputs and tests, together
with every declared variable, and put $s=|S|\ge1$. On $S$, draw an edge
from each immediate argument to its parent term. Vertices represent
expressions, not separate occurrences of those expressions. Let $\rho$
be the minimum size of a set meeting all directed paths from variables
to outputs, including paths of length zero. Repeated output terms
designate the same vertex.

Terms occurring only in tests do not change $\rho$: all ancestors of
an output already occur among its subterms. Constant symbols have no
incoming edges and start no path from a variable. Unused variables and
terms occurring only in tests can enlarge $S$ without increasing $\rho$.
For example, the outputs $(x,y,f(z),g(f(z)))$ have $\rho=3$, since the
last two terms share $f(z)$. The cut $\{x,y,f(z)\}$ has size three;
the paths consisting of $x$ alone, $y$ alone, and the edge from $z$ to
$f(z)$ show that no smaller cut is possible.

\begin{theorem}[The term-cut bounds with disequalities]\label{thm:cut}
If a test in $\Delta$ has syntactically identical sides, then
$\DD_n(\mathbf t,\Delta)=0$ for every $n\ge1$.
Otherwise, for every $n\ge s$,
\[
 \lfloor n/s\rfloor^\rho\le \DD_n(\mathbf t,\Delta)\le n^\rho.
\]
For each fixed instance without such a test,
\[
 \DD_n=\Theta(n^\rho),\qquad
 \lim_{n\to\infty}\frac{\log\DD_n}{\log n}=\rho.
\]
One can check for identical sides and compute $\rho$ in polynomial time
in the explicit input length.
\end{theorem}
\begin{proof}
A test $u\ne u$ cannot be satisfied. For the upper bound, fix an
interpretation and a cut $K$. The values of its terms determine every
output value. Indeed, evaluate an output recursively, stopping whenever
a vertex of $K$ is reached and using its given value. Any other leaf
must be a constant: reaching a variable without meeting $K$ would give
a path missed by the cut. Consequently two assignments with the same
values at $K$ have the same output tuple. There are at most $n^{|K|}$
possible tuples of values at $K$, and restricting the allowed inputs
cannot increase this bound. Take a minimum cut.

For the lower bound, Menger's theorem gives $\rho$ vertex-disjoint
paths from variables to outputs. Choose a set $B$ of size $m\ge1$ and
use the finite set $A=S\times B$. An element $(u,b)$ has a label
$u\in S$ and a value $b\in B$. We construct the functions so that
every subterm evaluates with its own label and the values at the
chosen path starts are copied along their paths.

For a symbol $f$ of arity $d$, consider arguments
$((u_1,b_1),\ldots,(u_d,b_d))$. If
$v=f(u_1,\ldots,u_d)$ belongs to $S$, let the result have label $v$.
If a chosen path enters $v$ through $u_j$, set its value to $b_j$;
otherwise use a fixed element of $B$. Vertex-disjointness makes the
chosen predecessor unique. If the predecessor term occurs in several
argument positions, select any one of those positions. On the assignments
used below they have equal values. For a constant symbol belonging to $S$, use the same
default value and its own label. All remaining function values may be
chosen arbitrarily. These instructions define one function for each
symbol: the symbol and the ordered input labels determine the term $v$,
so different subterms cannot impose conflicting instructions.

Give each source variable its own label. Allow the values at the
$\rho$ path starts to vary independently in $B$, and fix all other
source values. Induction on terms shows that each subterm has its own
label. Every test with different sides is therefore satisfied. Along
each selected path its starting value is copied to the output, giving
$m^\rho$ distinct output tuples. Hence $\DD_{sm}\ge m^\rho$.
For a larger set, extend each function arbitrarily while preserving
its values on $(S\times B)^d$. Taking $m=\lfloor n/s\rfloor$ proves
the lower bound for every $n\ge s$.

For $n\ge2s$, $\lfloor n/s\rfloor\ge n/(2s)$, giving the growth
and logarithmic limits, also when $\rho=0$. To compute the cut, split
each vertex into an entrance and exit joined by an edge of capacity one.
Replace each original edge by an edge from the corresponding exit to
entrance, of capacity $s+1$. Attach a new source to all variable
entrances and all output exits to a new sink, again with capacity
$s+1$. Cutting all unit edges costs $s$, so a minimum cut uses only
unit edges. Integral max-flow/min-cut gives both $\rho$ and the
vertex-disjoint paths in polynomial time. This is the directed vertex
form of Menger's theorem \cite{Menger}.
\end{proof}

\begin{proposition}[The sharper bound when $x,y$ are outputs]\label{prop:sharp-cut}
For every $\Gamma\in\mathcal C$ and $n\ge2$,
\[
 D_\Gamma(n)\le b_{\rho_\Gamma}(n)
                  =n^{\rho_\Gamma}-n^{\rho_\Gamma-1}.
\]
For every $j\ge2$ the bound is sharp: the tuple of $j$ distinct variables,
with only the test $x\ne y$ on its first two coordinates, has
$\rho_\Gamma=j$ and $D_\Gamma(n)=b_j(n)$.
\end{proposition}
\begin{proof}
Every cut $K$ contains $x,y$ because of their paths of length zero.
For a fixed interpretation, the values at $K$ determine the output
tuple, as in Theorem~\ref{thm:cut}. On allowed inputs the values at
$x,y$ differ, so there are at most
\[
 n(n-1)n^{|K|-2}=b_{|K|}(n)
\]
possible tuples of values at the cut. Choose a minimum cut and maximise
over interpretations. In the stated example, every allowed input is
a different output, and every one of the $j$ variables is an output.
Thus the cut has size $j$ and the image has size $b_j(n)$.
\end{proof}

\begin{corollary}[Thresholds for all sufficiently large sets]\label{cor:eventual}
Fix an integer $d\ge0$ and a function $h_d$ such that
$n^d<h_d(n)$ for all sufficiently large $n$, and
$h_d(n)=o(n^{d+1})$. Given an instance, one can decide in polynomial
time whether $\DD_n\ge h_d(n)$ for all sufficiently large $n$.
The answer is yes exactly when no test has syntactically identical
sides and $\rho\ge d+1$.
\end{corollary}
\begin{proof}
An identical-side test makes every image empty. Otherwise the lower
bound is at least a positive constant, depending on the instance,
times $n^\rho$ for all sufficiently large $n$. If $\rho\ge d+1$,
this eventually exceeds $h_d(n)$. If $\rho\le d$, the upper bound
$\DD_n\le n^d<h_d(n)$ gives the converse.
\end{proof}

\paragraph{How the labels handle repeated symbols.}
For terms $f(x,y)$ and $f(y,x)$ using the same binary symbol, the
ordered input labels are $(x,y)$ and $(y,x)$. They are different, so
one function can give each result its required label. Two occurrences
of the same term $f(x,y)$ instead use the same rule. For $f(u,u)$,
either argument position can be used to copy the value of $u$.
The rules depend on terms, never on occurrences independently.

This is the routing method of Riis and Gadouleau \cite{RG}, with the
additional observation that the labels satisfy all disequalities at
once. For example, $u\ne v$, $v\ne w$, and $w\ne u$ hold simultaneously
when $u,v,w$ are distinct terms and each receives its own label.
Small finite sets can still be too small for the tests. The theorem
controls growth as $n$ increases; it does not decide whether the sharp
upper bound is ever reached at a finite size.
\section{Forcing identities from all distinct pairs}

\begin{lemma}[An inequality characterises the identity function]\label{lem:avoid}
For any set $A$ and function $F:A\to A$,
\[
 [\forall a,b\in A\ (a\ne b\Rightarrow F(a)\ne b)]
 \quad\Longleftrightarrow\quad F=\id_A.
\]
\end{lemma}
\begin{proof}
If $F(a)\ne a$, take $b=F(a)$ to contradict the stated condition.
Conversely, if $F=\id_A$, the condition is just $a\ne b\Rightarrow a\ne b$.
\end{proof}

For functions $F_1,\ldots,F_M,W:A\to A$, define the set of allowed pairs
\[
 C(F,W)=\{(a,b)\in A^2:a\ne b,\ F_j(a)\ne b\ (1\le j\le M),\ W(a)\ne a\}.
\]
When $(x,y)$ are the only variables and are also the outputs, this set
is itself the output image. Its size is at most $n(n-1)$ for $|A|=n$.
The next proposition characterises equality, called saturation above.

\begin{proposition}[Equality with the pair bound]\label{prop:saturation}
If $|A|=n\ge2$, then
\[
 |C(F,W)|=n(n-1)
 \quad\Longleftrightarrow\quad
 [F_j=\id_A\text{ for all }j]\ \text{and}\ [W\text{ has no fixed point}].
\]
\end{proposition}
\begin{proof}
The set $C(F,W)$ is contained in the set of all $n(n-1)$ pairs with
different coordinates. Equality of their sizes means that every such
pair is allowed. Lemma~\ref{lem:avoid} then gives $F_j=\id_A$ for
each $j$. For any $a\in A$, choose $b\ne a$; the test on $W$ gives
$W(a)\ne a$. Conversely, these identities and the absence of fixed
points make every pair with different coordinates satisfy every test.
\end{proof}

We can also count exactly how many pairs are excluded. For each $a\in A$,
put $B_a=\{a,F_1(a),\ldots,F_M(a)\}$, and put
$Z=\{a\in A:W(a)=a\}$. If $a\in Z$, no pair with first coordinate
$a$ is allowed. Otherwise the allowed second coordinates are precisely
$A\setminus B_a$. Hence
\begin{align}
 |C(F,W)|&=\sum_{a\notin Z}(n-|B_a|),\label{eq:count}\\
 n(n-1)-|C(F,W)|&=(n-1)|Z|+
                  \sum_{a\notin Z}(|B_a|-1).\label{eq:defect}
\end{align}
If $W$ has no fixed point, then $Z$ is empty and $|B_a|-1\le M$, so
\begin{equation}\label{eq:densitybound}
 n(n-1)-Mn\le |C(F,W)|\le n(n-1).
\end{equation}
The lower bound is an ordinary integer inequality and remains valid
when $n(n-1)-Mn$ is negative.

\begin{example}[A failed identity can exclude just one pair]
Take $A=\{0,1,2\}$, outputs $(x,y)$, and tests $x\ne y$, $F(x)\ne y$.
If $F(0)=1$, $F(1)=1$, and $F(2)=2$, exactly five of the six pairs with
different coordinates are allowed: only $(0,1)$ is excluded. Admitting
all six pairs requires $F(0)=0$ as well.
If we also require $F(x)\ne x$, all six can never be admitted.
More generally, on any $n$-element set with $n\ge2$, these three tests
have maximum output image size $n(n-2)$. A point moved by $F$ allows
$n-2$ second coordinates, while a fixed point allows none. A cyclic
permutation of all $n$ elements attains $n(n-2)$.
Thus the proportion $n(n-2)/[n(n-1)]$ tends to one, although the image
always falls short of $n(n-1)$ by at least $n$.
\end{example}

\section{Reduction from a fixed finitely presented group}\label{sec:reduction}

\paragraph{The external undecidability theorem.}
A finite quotient of a group $G$ means a finite group $Q$ together with
a surjective homomorphism $G\to Q$. A word is \emph{separated} in such
a quotient if its image is not the identity. Bridson and Wilton
\cite[Theorem 2.1]{BW}, deriving a consequence of Slobodskoi's construction
\cite{Slo}, prove that there is a fixed finitely presented group $G$
for which the set
\[
 \{w:\text{$w$ has nonidentity image in some finite quotient of $G$}\}
\]
is recursively enumerable and undecidable. Thus an algorithm can
eventually find a separating finite quotient whenever one exists,
but no algorithm decides for every input word whether one exists.
The presentation is fixed; the input consists only of the word.
We first describe our reduction for an arbitrary finite presentation.

Fix
\[
 G=\langle a_1,\ldots,a_d\mid r_1,\ldots,r_q\rangle.
\]
For each generator introduce two unary symbols $p_i,m_i$. The notation
suggests a generator and its inverse, but their interpretations are
initially independent functions. For a signed word
$v=\ell_1\cdots\ell_L$, whose letters belong to
$\{a_1^{\pm1},\ldots,a_d^{\pm1}\}$, define
\[
 \tau_v(x)=f_{\ell_1}(f_{\ell_2}(\cdots f_{\ell_L}(x)\cdots)),
 \qquad \tau_\varepsilon(x)=x,
\]
where $f_{a_i}=p_i$, $f_{a_i^{-1}}=m_i$, and $\varepsilon$ is the
empty word. Let $\mathcal L_G$ be the list of $M=2d+q$ terms
\[
 p_i(m_i(x)),\ m_i(p_i(x))\quad(1\le i\le d),
 \qquad \tau_{r_j}(x)\quad(1\le j\le q).
\]
We call these the \emph{law terms}, since requiring each to equal $x$
would impose inverse and relator identities. Repetitions are retained;
an empty relator contributes the term $x$.
Our convention for multiplying permutations is composition with the
rightmost function applied first, agreeing with the displayed nesting.

Given the input word $w$, construct an instance with just the variables
$x,y$, output tuple $(x,y)$, and tests
\begin{equation}\label{eq:compiler}
 x\ne y,\qquad L(x)\ne y\quad(L\in\mathcal L_G),\qquad
 \tau_w(x)\ne x.
\end{equation}
Here $L(x)$ denotes the law term $L$, with its indicated variable $x$.
Write $D_{G,w}(n)$ for this instance's dispersion.
For each interpretation its output image is $C(F,W)$ above, where
$F_1,\ldots,F_M$ interpret the law terms and $W$ interprets $\tau_w$.

\begin{theorem}[The reduction]\label{thm:reduction}
For every finite presentation $G$ and signed word $w$,
\[
\begin{gathered}
 \exists n\ge2:\ D_{G,w}(n)=n(n-1)\\
 \Longleftrightarrow\\
 \text{$w$ has nonidentity image in some finite quotient of $G$}.
\end{gathered}
\]
The construction has linear output length in the explicit presentation
and word lengths, and is computable in polynomial time.
\end{theorem}
\begin{proof}
Suppose $D_{G,w}(n)=n(n-1)$, and choose an interpretation attaining
this maximum on a set $A$ of size $n$. Write $p_i,m_i$ also for its
unary functions. Proposition~\ref{prop:saturation} gives
$p_i\circ m_i=\id_A=m_i\circ p_i$, so the $p_i$ are permutations and
$m_i=p_i^{-1}$. It also gives $\tau_{r_j}=\id_A$ for every relator.
The universal property of the presentation therefore gives a
homomorphism $\phi:G\to\Sym(A)$. The interpreted $\tau_w=\phi(w)$
has no fixed point, hence is not the identity. The finite group
$\phi(G)$, with the surjection $G\to\phi(G)$, is the required quotient.
The identities hold at every point of $A$, as required to define this
homomorphism.

Conversely, suppose $\phi:G\to Q$ is a finite quotient with
$\phi(w)\ne1$. On $A=Q$, interpret $p_i$ as left multiplication by
$\phi(a_i)$ and $m_i$ as left multiplication by its inverse.
Each law term is the identity function. The term $\tau_w$ acts by
left multiplication by $\phi(w)$, which has no fixed point:
$\phi(w)b=b$ would imply $\phi(w)=1$. Proposition~\ref{prop:saturation}
gives $D_{G,w}(|Q|)=|Q|(|Q|-1)$, and $|Q|\ge2$.
The construction of the instance only lists the inverse-law terms,
copies the relators, and writes the nested term for $w$, giving the
stated size and time bounds.
\end{proof}

\begin{theorem}[Undecidability with fixed unary symbols and tests]
\label{thm:fixed-signature}
There is a fixed finite presentation $G$ such that, given a signed
word $w$, deciding whether $D_{G,w}(n)=n(n-1)$ for some $n\ge2$ is
undecidable. The associated instances use a fixed finite unary
signature, a fixed number of tests, and just two variables and two
outputs. Their positive instances are recursively enumerable but not
decidable, and their complement is not recursively enumerable.
\end{theorem}
\begin{proof}
Fix the presentation provided by the external theorem above and apply
Theorem~\ref{thm:reduction}. Only $\tau_w$ varies with the input.
For enumerability, try $n=2,3,\ldots$ in order, and for each $n$
enumerate all interpretations on $[n]$. There are finitely many at
each size, and each output image is computable by finite enumeration.
This procedure finds an image of size $n(n-1)$ whenever one exists.
If the negative instances also had such a semidecision procedure,
running both procedures in parallel would decide every instance,
contradicting the reduction.
\end{proof}

Bridson--Wilton state their theorem for freely reduced words, meaning
words with no adjacent inverse pair $a_i a_i^{-1}$ or $a_i^{-1}a_i$.
Such a word is already a valid signed input word here. Conversely,
removing adjacent inverse pairs is computable and preserves evaluation
in every group, hence preserves separation in finite quotients.
The reduction theorem therefore gives the same exact decision problem
for arbitrary signed words.
This does \emph{not} license cancellation in the terms under arbitrary
interpretations: $p_i$ and $m_i$ need not be inverse functions.
Cancellation preserves their evaluation only after the required inverse
laws hold, as they do when the pair bound is attained.

The empty word yields the impossible test $x\ne x$ and is a decidable
negative case. Restricting to nonempty freely reduced words consequently
leaves the problem undecidable. Recursively enumerable completeness
is not asserted; it does not follow from enumerability and
undecidability alone.

\section{Proof of the one-unit theorem}\label{sec:mainproof}
\begin{proof}[Proof of Theorem~\ref{thm:main}]
Fix $k\ge2$. Start with the instance associated to $G,w$ in
Theorem~\ref{thm:reduction}. Add $k-2$ new variables as further output
coordinates, with no tests involving them. This operation is called
\emph{padding}. For each interpretation the new image is exactly the
Cartesian product of the old pair image with $[n]^{k-2}$: the new
variables can be chosen independently of every test and of the old
outputs. Maximising over the same interpretations therefore gives
\[
 D_{G,w}^{(k)}(n)=n^{k-2}D_{G,w}(n),
\]
where the superscript denotes the padded instance.
All $k$ variables are outputs, so every cut contains them and they
also form a cut. Thus $\rho=k$.
Since $D_{G,w}(n)\le n(n-1)$, for $n\ge2$ the inequality
$D_{G,w}^{(k)}(n)\ge b_k(n)$ holds exactly when
$D_{G,w}(n)=n(n-1)$. For fixed $k$, adding the variables and outputs
takes constant extra space; the signature and tests are unchanged.
Theorem~\ref{thm:fixed-signature} proves the claimed undecidability.

For the strict problem, a test with identical sides makes the answer
no. Otherwise compute $\rho=\rho_\Gamma$. If $\rho\le k$, then
Proposition~\ref{prop:sharp-cut} gives, for every $n\ge2$,
\[
 D_\Gamma(n)\le b_\rho(n)\le b_k(n),
\]
so again the answer is no. If $\rho\ge k+1$, Theorem~\ref{thm:cut} gives
\[
 D_\Gamma(n)\ge\lfloor n/s\rfloor^\rho
                 \ge\lfloor n/s\rfloor^{k+1}\qquad(n\ge s).
\]
For $n=s^{k+1}$ the last expression is
\[
 (s^k)^{k+1}=n^k\ge n^k-n^{k-1}+1=b_k(n)+1.
\]
Here $s\ge2$ because $x,y$ are distinct, so $n\ge2$.
This proves the criterion and the explicit sufficient size.
Checking syntax and computing a minimum vertex cut are polynomial-time
operations, as detailed below. Comparing $\rho$ with the binary
integer $k+1$ is polynomial too. The decision algorithm need not
construct $n$ or any interpreting functions.

Finally, suppose a total computable bound $B(N)$ as in the theorem
existed for some fixed $k$. On an input from the padded family,
compute its bit length $N$ and enumerate $2\le n\le B(N)$.
For each size enumerate all interpretations and all assignments,
form the set of distinct outputs satisfying the tests, and compare its
size with $b_k(n)$. The search is finite and, by the assumed bound,
would answer yes exactly on the positive instances. This would decide
the undecidable problem. The contradiction needs no time bound on
this exhaustive search.
\end{proof}

\paragraph{The polynomial-time algorithm.}
First validate the input as described in Section~\ref{sec:language}.
Process terms in topological order. Assign an identifier to each source
and identify compound terms by the record $(f,i_1,\ldots,i_d)$ of their
symbol and ordered child identifiers. Equal records receive the same
identifier. Thus syntactically equal expressions are identified even
when the input supplied different node names for them.
Reject any test whose two resulting identifiers agree. Build the term
graph on these identifiers, compute its minimum cut by the network in
Theorem~\ref{thm:cut}, and answer yes exactly when $\rho_\Gamma\ge k+1$.

For a term description of $N$ bits, at most $N$ nodes and argument
entries are supplied. Pairwise comparison of the records can be done
in $O(N^3)$ bit operations. The split network has $O(N)$ vertices and
edges and integer capacities at most $N+1$. An augmenting-path flow
algorithm needs at most $s\le N$ augmentations, since each increases
the integer flow by at least one, while the cut of all unit edges
bounds it by $s$. Path searches and updates involve only polynomially
many entries. Indices, capacities and flow values have $O(\log N)$
bits, so even sequential scans for lookups give polynomial bit cost.
Together with validation and the binary comparison with $k+1$, this
proves the uniform polynomial-time bound. No loop of length $k$ is
needed when $k$ is large.

\paragraph{Which degree can give an undecidable question.}
For an instance in $\mathcal C$ without an identical-side test,
$D_\Gamma(n)=\Theta(n^\rho)$. If $k<\rho$, the strict threshold is
exceeded for all sufficiently large $n$. If $k>\rho$, the sharp upper
bound is smaller than $b_k(n)$ for every $n\ge2$, so even the non-strict
answer is no. When $k=\rho$, the non-strict question asks whether the
upper bound $b_\rho(n)$ is attained for some $n$; this is the subclass
on which the reduction proves undecidability.
It also follows that the existential strict question is equivalent
to requiring $D_\Gamma(n)\ge b_k(n)+1$ for all sufficiently large $n$,
and to requiring $D_\Gamma(n)\ge n^{k+1/2}$ for all sufficiently large
$n$. These are consequences; both questions in Theorem~\ref{thm:main}
itself use the existential quantifier.

\paragraph{Three outputs.}
For outputs $(x,y,t)$ in $\mathcal C$, first reject a test with identical
sides. Then $\rho=2$ if $t$ contains no variable other than $x,y$,
including when it contains no variables at all. If another variable
occurs in $t$, a path from it to $t$ cannot pass through $x$ or $y$,
since variable vertices have no incoming edges. Together with the
paths of length zero at $x,y$, this gives three disjoint paths.
The three outputs form a cut, so $\rho=3$.
For $k=2$, the strict algorithm thus only has to check whether $t$
contains another variable. In the positive case $s^3$ is a sufficient
size and $D_\Gamma(n)=\Theta(n^3)$. The lower bound is
$n^3/(8s^3)$ for $n\ge2s$, hence at least $n^{5/2}$ for $n\ge64s^6$.

Requiring $x,y$ to be outputs is essential to this criterion.
For example, with four distinct variables $x,y,z,w$, outputs $(z,w)$,
and only the test $x\ne y$, every pair $(z,w)$ is an output for
$n\ge2$. Thus $D(n)=n^2\ge b_2(n)+1$ although the cut size is two.
This instance lies outside $\mathcal C$: the two variables constrained
to differ are not output coordinates.

\section{Functions approaching the pair bound}\label{sec:near}

\begin{theorem}[An explicit approximation to the pair bound]\label{thm:near}
Fix a finite presentation $G$ and its $M=2d+q$ law terms.
For every nonempty signed word $w$ of length $L$ and every integer $n>L$,
\[
 n(n-1)-Mn\le D_{G,w}(n)\le n(n-1).
\]
An interpretation giving the lower bound is explicitly computable and
uses only cyclic shifts. In particular, for each such word,
\[
 \lim_{n\to\infty}\frac{D_{G,w}(n)}{n(n-1)}=1,
 \qquad D_{G,w}(n)=n^2-O_G(n).
\]
The constant in the error bound depends only on the presentation,
while the condition $n>L$ also depends on the word.
\end{theorem}
\begin{proof}
On $A=\mathbb Z/n\mathbb Z$, interpret every $p_i$ and every $m_i$ as
$z\mapsto z+1$. A word of length $\ell$ as written then acts by
$z\mapsto z+\ell$, regardless of the signs of its letters.
Because $0<L<n$, the function interpreting $\tau_w$ has no fixed point.
For each first coordinate $a$, the tests exclude at most $M+1$ second
coordinates: $a$ itself and the $M$ law values. The lower bound follows
from \eqref{eq:densitybound}. Dividing by $n(n-1)$ gives a proportion
between $1-M/(n-1)$ and $1$, proving the limit and the error estimate.
The inequality is valid even if its lower bound is nonpositive.

The formal inverse symbols receive the same shift as the positive
symbols. This is permitted because an interpretation consists of
arbitrary functions. Inverse laws are forced only when every pair with
different coordinates is admitted.
\end{proof}

For this interpretation the exact count is
\[
 |C|=n(n-1-|R_n|),\qquad
 R_n=\{|v|\bmod n:\tau_v\in\mathcal L_G\}\setminus\{0\}.
\]
Each inverse-law term is here regarded as its two-letter word, and
$|v|$ always counts the letters without cancellation. Distinct nonzero
residues give distinct additional forbidden second coordinates for
each first coordinate; the zero residue forbids nothing beyond $a$.

Given a nonempty word $w$ and a positive rational $\varepsilon$, choose,
for example,
\[
 n=1+\max\{L,\lceil1+M/\varepsilon\rceil\}.
\]
This integer satisfies $n>L$ and $M/(n-1)\le\varepsilon$.
Return $n$ and the unary functions $z\mapsto z+1\pmod n$.
Under this interpretation the instance has at least
$(1-\varepsilon)n(n-1)$ distinct output pairs.
The functions have finite explicit tables; writing those tables takes
time polynomial in $n$ and the input length, with $n$ measured by its
value rather than its binary length. Here \emph{density} refers exactly
to the ratio of the image size to $n(n-1)$.

This gives more information than the quadratic growth exponent:
for any prescribed rational error an algorithm produces functions with
the required proportion of outputs. Nonetheless it cannot decide
whether some interpretation admits all $n(n-1)$ pairs. A nonzero loss
of at most $Mn$ pairs becomes arbitrarily small after division by
$n(n-1)$, while zero loss imposes every group identity on all of $A$.
The approximate interpretations are not required to be group actions.

The distinction between an empty word and a word representing the
identity matters. The empty word gives $x\ne x$ and zero dispersion.
A nonempty word can represent the identity in $G$, or even cancel
freely to the empty word, while retaining positive length as an input
expression. Its shift interpretation still gives the bound above.
There is no conflict with the reduction theorem: cancellation preserves
group evaluation, but need not preserve evaluation by these arbitrary
unary functions. For such an identity word, the exact answer is no
even though the proportions tend to one.

For the padded instances of degree $k\ge2$, the exact product identity
gives, for $n>L$,
\[
 b_k(n)-Mn^{k-1}\le D_{G,w}^{(k)}(n)\le b_k(n),
 \qquad
 \lim_{n\to\infty}\frac{D_{G,w}^{(k)}(n)}{b_k(n)}=1.
\]
The lower bound again uses ordinary integer subtraction. For every
$n\ge2$ the ratio equals the original pair ratio. Thus the same
algorithm gives a proportion at least $1-\varepsilon$ of $b_k(n)$
by adding the independent output variables, although deciding whether
$b_k(n)$ is attained at any size remains undecidable.

\section{Scope and the problem without tests}\label{sec:open}
The decision problems use one finite set for all variables and function values.
For each fixed $k\ge2$, undecidability already holds with $k$ variables,
all retained as outputs, a fixed finite unary signature and a fixed
number of tests. Only the nested expression for the input word varies.
The tests compare computed expressions as well as variables.
Writing a computed term as a new variable with a defining equation does
not make it an independent source variable: its value is determined
by the original inputs. The reduction therefore proves no
undecidability result for the narrower language in which tests may
compare only source variables.

\paragraph{The open bijectivity problem.}
Without disequality tests, finite square term bijectivity asks the
following: given $k$ terms in $k$ variables, does some interpretation
on a finite set $A$, with $|A|\ge2$, make the resulting map
$A^k\to A^k$ bijective? This problem remains open in general.
The analysis of ordinary dispersion in \cite{Paper2} identifies an image
of size $|A|^k$ for such a map with bijectivity, and gives padding
reductions from the case $k=3$
to the problem of attaining the full output set at every fixed output
dimension at least three. Those reductions do not prove undecidability
and are not used here. An unconditional polynomial-time/undecidable
separation for dispersion without disequality tests remains open.
The present bounds also concern a different quantity from the number
of solutions of an arbitrary equation system or the entropy and
guessing number of a graph; they do not imply that those quantities
are integer-valued or efficiently computable.

\paragraph{Machine verification.}
The development, archived release and interactive paper are available at:
\begin{quote}
\noindent\textbf{GitHub:}\\
\PaperThreeLeanGitHub\\
\textbf{DOI:} \PaperThreeLeanDOI\\
\textbf{Interactive paper:}\\
\PaperThreeInteractive
\end{quote}
The accompanying Lean 4.26 development has 71 mathematical and
algorithmic modules. It checks the output-image definition, the sharp
cut bound when $x,y$ are outputs, finite directed Menger and max-flow
results, the routing construction, and the strict criterion for every
degree with its explicit sufficient size. It also checks the argument
forcing identities, the equivalence with finite-quotient separation,
the exact product under padding, finite search, enumerability, and the
absence of a computable bound on alphabet size, as well as the
approximation, asymptotic and flattening results.
The undecidability theorems assume explicitly the fixed-presentation
finite-group noncomputability result. The classical theorem of
Slobodskoi in the form used by Bridson--Wilton is an external input,
not a theorem formalised in this project. The axiom audits use only
the standard axioms of Lean's classical foundations, with no admissions.
The build and audit records accompany the source.

For the strict decision algorithm, the development proves correctness
of parsing, rejection of malformed bit strings, and the decision by
max-flow. It also certifies one fixed deterministic random-access
machine program under a logarithmic cost criterion. For every input
bit string $b$, including malformed inputs, this program halts with
the specified decision result. Fixed constants $C\ge1$ and $c$ bound
both its number of instructions executed and its total logarithmic
cost by $C(|b|+1)^c$. The binary degree is part of the input; the
program does not loop to its value. The main formal statements are
\texttt{program\_polynomial} and \texttt{program\_accepts\_iff}
in \texttt{MMain.lean}. The standard polynomial simulation from
this RAM model to a multitape Turing machine is not formalised here.
The ordinary bit-complexity argument appears in
Section~\ref{sec:mainproof}.

The interactive edition links the text of this paper to explanations, exact formal
statements and proofs, including the hypotheses used at each stage.
The reader can inspect Lean's recorded proof states and follow formal
dependencies into library proofs and the underlying logical rules.
External mathematical input and the limits of machine verification
are identified explicitly. This correspondence helps a reader check
that the formal statements express the claims made in the paper;
kernel checking of a Lean proof alone does not establish that
correspondence.

\bibliographystyle{plain}
\bibliography{references}
\end{document}